\documentclass[letterpaper]{article}

\PassOptionsToPackage{numbers}{natbib}

\usepackage[preprint]{neurips_2025}

\usepackage[utf8]{inputenc}
\usepackage[T1]{fontenc}
\usepackage{hyperref}
\usepackage{url}
\usepackage{booktabs}
\usepackage{amsfonts}
\usepackage{nicefrac}
\usepackage{microtype}
\usepackage{xcolor}
\usepackage{threeparttable}
\usepackage{tabularx}

\usepackage{amsmath, amssymb, mathtools, amsthm, mathrsfs, bbm}
\usepackage{booktabs, multirow, siunitx, pdflscape, float}
\usepackage{microtype, subcaption, multirow, makecell, calc, comment, url}
\usepackage[pdftex]{graphicx}
\usepackage{algorithm, algpseudocode}
\newfloat{implementation}{htbp}{lop}
\floatname{implementation}{Implementation}

\newtheorem{theorem}{Theorem}

\newtheorem{example}{Example}

\newcommand{\best}[1]{\multicolumn{1}{c}{\textbf{\boldmath $#1$}}}

\newcommand\cE{\mathcal{E}}
\newcommand\cF{\mathcal{F}}
\newcommand\cG{\mathcal{G}}

\newcommand\bR{\mathbb{R}}

\title{Fractional-Boundary-Regularized Deep Galerkin Method for Variational Inequalities in Mixed Optimal Stopping and Control}

\author{
  Yun~Zhao \\
  Department of Mathematics\\
  Imperial College London\\
  London, SW7 2AZ \\
  \texttt{yun.zhao23@imperial.ac.uk} \\
  \And
  Harry Zheng \\
  Department of Mathematics\\
  Imperial College London \\
  London, SW7 2AZ \\
  \texttt{h.zheng@imperial.ac.uk} \\
}

\begin{document}

\maketitle

\begin{abstract}
    Mixed optimal stopping and stochastic control problems define variational inequalities with non-linear Hamilton-Jacobi-Bellman (HJB) operators, whose numerical solution is notoriously difficult and lack of reliable benchmarks. We first use the dual approach to transform it into a linear operator, and then introduce a Fractional-Boundary-Regularized Deep Galerkin Method (FBR-DGM) that augments the classical $L^2$ loss with Sobolev-Slobodeckij norms on the parabolic boundary, enforcing regularity and yielding consistent improvements in the network approximation and its derivatives. The improved accuracy allows the network to be converted back to the original solution using the dual transform. The self-consistency and stability of the network can be tested by checking the primal-dual relationship among optimal value, optimal wealth, and optimal control, offering innovative benchmarks in the absence of analytical solutions.
\end{abstract}

\section{Introduction}
\textbf{Limitations of existing methods.} Optimal stopping mixed with stochastic control problems are a class of challenging problems in quantitative finance, operations research and reinforcement learning. Mathematically, these problems lead to a non-linear variational inequality which does not have a closed-form solution in general. One way to solve these problems is to use the dual approach, which transforms the non-linear variational inequality into a linear one. However, the complexity of the linear operator remains very high, as evidenced by the fact that the corresponding free boundaries can have multiple shapes and occur in varying numbers, see Figure 1 in \citet{maGlobalClosedFormApproximation2019}. Another aspect of the dual approach is that the derivatives of the dual value function has to converge in some way to make sure the primal-dual transform does not ruin the convergence of the original problem. This is a very tough requirement, as common numerical algorithms for variational inequalities only give the solution based on a stationary initial condition, and the derivatives can only be approximated by Finite Difference Method (FDM) with many repetitions for small changes, e.g. \citet{beckerDeepOptimalStopping2019}, \citet{beckerSolvingHighdimensionalOptimal2021}, \citet{longstaffValuingAmericanOptions2001}. Empirical evidence shows that these approaches generally give very poor approximations of the derivatives, and thus the primal-dual transform is not reliable.

\textbf{Motivations for using deep learning.} Instead, we propose to use deep learning to solve the mixed optimal stopping and control problem. The main idea is to use a neural network to approximate the solution to the dual problem, and then transform it back to the original problem using the dual transform. The advantage of this approach is that the trained network is a function of the whole domain, and the derivatives at any point can be computed directly and instantly. This is much more efficient than repeating the traditional algorithms and using FDM to get derivatives at only one point. Moreover, we do not need to consider the properties of free boundaries, continuation region and stopping region, as we are approximating the solution using the minimum representation of the variational inequality, which holds on the whole domain.

\textbf{Main challenges.} The only problem is that the convergence of the network is not guaranteed. The well-known universal approximation theorem for deep neural networks only states the existence of a network that can approximate any $C^k$ function and its derivatives on a compact domain, see \citet{hornikUniversalApproximationUnknown1990}. Based on it, the Deep Galerkin Method (DGM) was proposed to solve PDEs, see \citet{sirignanoDGMDeepLearning2018}. Using loss functions that only include $L^2$ norms, DGM is a very general and accurate way to solve quasi-linear PDEs, and the convergence is proved under additional conditions. However, variational inequalities are very different from PDEs, as the operator part is only equal to zero on the continuation region, which is separated from the stopping region by an unknown free boundary. \citet{sirignanoDGMDeepLearning2018} approximates the free boundary iteratively by using the approximated solution. Strictly speaking, this approach is different from solving a PDE, as the convergence of the approximated free boundary needs to be examined, which is not done in the original paper. Moreover, it is less efficient than directly using the minimum representation in the loss function.

\textbf{Proposed solution.} In this paper, we propose a new loss function that augments the original $L^2$ norm with Sobolev-Slobodeckij norms on the parabolic boundary, also called fractional Sobolev norms. The idea is to match the network output with the boundary condition of the dual value function to an extent more than $L^2$, which can improve the interior approximations for higher derivatives. This is due to the inverse trace mapping from functional analysis, which states there is a continuous mapping from the fractional Sobolev space on boundary to the Sobolev space on whole domain with a higher regularity. The continuity of the mapping ensures that once the boundary norm converges to zero, the interior norm also converges to zero. Applying this to our case, we can improve the convergence of the derivatives of the network output, which is crucial for the primal-dual transform. The new loss function is called the Fractional-Boundary-Regularized loss, and we denote the network trained with it as FBR-DGM.

\textbf{Empirical validation.} With neural networks, we numerically test the approximations compared to the Binomial Tree Method (BTM), which is the most widely used algorithm in optimal stopping problems, and the Global Closed Approximation (GCA) method, which is a classical method proposed by \citet{maGlobalClosedFormApproximation2019} to solve the dual problem by approximating the free boundary. We also take the full advantage of the primal-dual relationship by checking the self-consistency of optimal value, optimal wealth and optimal control, all of which are computed from the network in some way. The optimal value is the output of the neural network transformed back to the original problem. The optimal wealth and optimal control are the derivatives of the network, which are computed directly from automatic differentiation. The numerical results show that, using BTM as the benchmark, both DGM and FBR-DGM are more accurate and efficient than GCA, and FBR-DGM further outperforms DGM in terms of accuracy and stability. The self-consistency check shows that it is easier to be consistent with power utility functions than with non-HARA utility functions, and FBR-DGM is more self-consistent than DGM in most cases.

\textbf{Main contributions.} We summarize our main contributions as follows:
\begin{enumerate}
\item We propose a new loss function that includes Sobolev-Slobodeckij norms on the parabolic boundary, which improves the convergence of the neural network and its derivatives. It is easy to implement in practice and shows powerful approximation capability in solving mixed optimal stopping and control problems.
\item We provide two innovative self-consistency checks using the primal-dual relationship, which can be used as extra tests for the accuracy of neural network output and are particularly important in the absence of analytical solutions.
\end{enumerate}

\section{Related work} \label{sec:related work}
\textbf{Deep learning for PDEs.} Many deep learning approaches have been proposed to solve PDEs, including the Deep Galerkin Method (DGM), \citet{sirignanoDGMDeepLearning2018} and physics-informed neural networks (PINNs), \citet{shinConvergencePhysicsInformed2020}. These methods achieve high-dimensional accuracy but consider the operator equals to be zero everywhere and do not solve the free boundary that defines the variational inequality.

\textbf{Deep learning for HJB equations.} One can solve HJB equations of continuous-time by feedforward neural networks. In \citet{nakamura-zimmererAdaptiveDeepLearning2021}, adaptive data-generation techniques extend these HJB nets to 30-plus dimensions and give real-time feedback laws. However, they assume actions are continuous, and as soon as a stopping decision is introduced, the HJB turns into a variational inequality and there is no means of capturing the free boundary.

\textbf{Deep learning for variational inequalities.} Little work exists that directly tackles variational inequalities via learning so far. \citet{hureDeepBackwardSchemes2020} tackles the non-linear variational inequality using reflected backward differential equations (RBSDE). Their convergence proof requires the Lipschitz condition on the driver and obstacle, and did not analyze the estimate error from neural network, which serves as an unknown term in the analysis. \citet{wangDeepLearningFree2021a} focuses on the Stefan problem and learns both the temperature field and the melting surface. However, it is only in 1-2 dimensional heat-equation settings, and it is needed to know the shape and numbers of the free boundary in advance. A very recent article \citet{alphonseNeuralNetworkApproach2024} casts the classical elliptic obstacle problem as a min-max representation between two neural networks and derives a priori error estimates. While their weak-adversarial formulation is elegant and can handle non-symmetric elliptic operators, it remains restricted to stationary (time-independent) operators, and there is no parabolic or HJB structure, so their approach cannot handle mixed optimal stopping and control problems evolving over time.

\textbf{Non-learning algorithms for HJB variational inequalities.} \citet{maGlobalClosedFormApproximation2019} analytically tackles the free boundary for a specific mixed stop and control model, but does not give the function on the whole domain, and requires FDM in order to compute the derivatives. Penalty plus policy-iteration schemes \citet{reisingerPenaltySchemePolicy2021} provides convergent discretization for HJB variational inequalities with jumps. It uses discretized schemes to approximate the penalized solutions, which introduce the corresponding discretization errors. Computing gradients or Hessians elsewhere requires additional finite-difference or re-solving on refined grids, which is costly and can amplify numerical noise. The convergence of derivatives is not examined.

\begin{table}[htbp]
    \centering
    \caption{Position of related works in the mixed stop and control problems}
    \label{tab:related_work}
    \resizebox{\linewidth}{!}{
    \begin{tabular}{@{} l c l >{\raggedright\arraybackslash\sloppy}p{4cm} >{\raggedright\arraybackslash\sloppy}p{4cm} @{}}
      \toprule
      Category &
      \makecell{Representative\\paper} &
      Operator type &
      Main Contributions &
      Limitation for mixed stop and control \\[2pt]
      \midrule

      Deep PDE &
      \citet{sirignanoDGMDeepLearning2018} &
      Parabolic &
      $L^2$ loss; Convergence for quasi-linear PDEs &
      Not for variational inequalities \\[4pt]
      
      Deep HJB equations &
      \citet{nakamura-zimmererAdaptiveDeepLearning2021} &
      Parabolic HJB &
      Adaptive DL; Nonlinear &
      Not for variational inequalities \\[4pt]

      Deep variational inequalities &
      \citet{hureDeepBackwardSchemes2020} &
      Nonlinear Parabolic &
      Two nets for value and gradient; RBSDE &Lipschitz; No error analysis for the neural network \\[4pt]
      
      Deep variational inequalities &
      \citet{wangDeepLearningFree2021a} &
      Heat equation &
      Two nets for value and free boundary; Stefan interface &
      Heat-equation specific; Prior knowledge to the free boundary; No control \\[4pt]

      Deep variational inequalities &
      \citet{alphonseNeuralNetworkApproach2024} &
      Elliptic &
      Two-net min-max gap formulation &
      Stationary (no time); No control \\[4pt]
      
      HJB variational inequalities &
      \citet{maGlobalClosedFormApproximation2019} &
      Parabolic HJB &
      Implicit free boundary for dual problem &
      Not a function on the whole domain \\[2pt]

      HJB variational inequalities &
      \citet{reisingerPenaltySchemePolicy2021} &
      Parabolic HJB &
      Convergent discretization; Infinite jumps &
      Pre-determined mesh; Convergence of derivatives not examined \\[4pt]

      \bottomrule
    \end{tabular}}
  \end{table}

\section{Preliminaries - mixed optimal stopping and control problem} \label{sec:preliminaries}
In this paper, we work on the particular example from \citet{maGlobalClosedFormApproximation2019}. The settings are briefly described as follows.

\textbf{Market model.} There is a complete market with a probability space $(\Omega, \mathcal{F}, P)$ and a natural filtration $(\mathcal{F}_t)$ generated by a standard Brownian motion $W$ which satisfies the usual conditions. There is one riskless savings account with interest rate $r > 0$ and one risky asset, namely, a stock. The stock process is $dS_t = \mu S_t \, dt + \sigma S_t \, dW_t,$
where $\mu, \sigma >0$ are the stock growth rate and volatility respectively. The wealth process is $dX_t = r X_t \, dt + \sigma \pi_t \left( \theta \, dt + dW_t \right),$
where $r>0$ is the riskless interest rate, $\theta = \frac{\mu - r}{\sigma}$ is the market price of risk. $(\pi_t)_{0 \leq t \leq T}$ denotes the portfolio process, and it is $\mathcal{F}_t$-progressively measurable satisfying $\mathbb{E}\left[\int_0^T |\pi_t|^2 dt\right] < \infty$.

\textbf{Primal value function.} For a minimum wealth threshold value $K>0$, a terminal investing time $T>0$, a utility discount factor $\beta > 0$, and an $\cF_t$-adapted stopping time $\tau \in [0,T]$, the primal value function is $V(t,x):(0,T)\times (K,+\infty) \to [0,+\infty)$, given by
\[
V(t, x) = \sup_{\tau, \pi} \mathbb{E} \left[ e^{-\beta(\tau - t)} U(X_\tau^{t, x, \pi} - K) \mid X_t = x \right], \quad (t,x) \in Q_x:=(0,T) \times (K,+\infty),
\]
where $U(x):[0,+\infty)\to [0,+\infty)$ is a utility function that satisfies the following conditions: $U \in C^2([0,\infty)),$ $U \text{ is increasing and strictly concave,}$ $U(0) = 0, \, \lim_{x\to\infty} U(x) = +\infty,$ $U'(0) = +\infty, \, \lim_{x\to\infty} U'(x) = 0,$ $U(x) < C(1+x^p) \text{ for } x \geq 0 \text{ and some } C>0, \, 0<p<1,$ $U(x) = -\infty \text{ for } x < 0.$
The primal value function $V(t,x)$ satisfies the following variational inequality:
\begin{equation} \label{eq: VI of V with control}
    \min \left\{-V_t - \sup_{\pi} \{rxV_x-\beta V+\pi(\mu-r)V_x+\frac{1}{2}\pi^2\sigma^2 V_{xx}\}, V-U(x-K)\right\}=0, \, (t,x) \in Q_x,
\end{equation}

Supposing that $V(t,\cdot)$ is strictly concave, solving the supremum over $\pi$, $V$ satisfies a non-linear HJB variational equation:
\begin{align} \label{eq: VI of V without control}
    \begin{cases}
        \min\{-V_t + \frac{\theta^2}{2} \frac{V_x^2}{V_{xx}} - rx V_x + \beta V, \, V - U(x-K)\} = 0, & \text{for } (t,x) \in Q_x, \\
        V(t, K) = 0, & \text{for } t \in (0,T), \\
        V(T,x) = U(x-K), & \text{for } x \in (K,+\infty), \\
    \end{cases}
\end{align}
It is difficult to solve \eqref{eq: VI of V without control}. Hence, we turn to the dual method described as follows.

\textbf{Dual value function.} Define the dual process $(Y_t)_{0\leq t \leq T}$ that satisfies
\begin{equation} \label{eq: dual process}
dY_t = (\beta - r) Y_t \, dt - \theta Y_t \, dW_t.
\end{equation}
The dual function of $U(\cdot - K)$ is given by
\begin{equation} \label{def: dual utility function}
\widetilde{U}_K(y) := \sup_{x>K}\{U(x-K)-xy\}=\widetilde{U}_0(y)-Ky \quad \text{ for } y>0.
\end{equation}
$\widetilde{U}_K$ is continuously differentiable, decreasing, strictly convex and satisfies $\widetilde{U}_K(0) = \infty, -Ky \leq \widetilde{U}_K(y) \leq \widetilde{C}y^{p-1} - Ky,$
where $\widetilde{C} = \max\{C, (Cp)^{\frac{1}{p-1}} [p^{-1}-1]\}$ and $C$, $p$ are from the conditions on $U$. The dual value function is defined as
\[
\widetilde{V}(t, y) := \sup_{t \leq \tau \leq T} \mathbb{E} \left[ e^{-\beta(\tau - t)} \widetilde{U}_K(Y_\tau) \mid Y_t = y \right], \quad (t,y) \in Q_y := (0,T)\times(0,+\infty). \]
We denote the continuation and stopping region as
\[
    C_y := \left\{(t,y) \in Q_y: \widetilde{V}(t,y) > \widetilde{U}_K(y)\right\}, \quad
    S_y := \left\{(t,y) \in Q_y: \widetilde{V}(t,y) = \widetilde{U}_K(y)\right\}.
\]
From dynamic programming principle, $\widetilde{V}(t,y)$ satisfies a linear variational inequality:
\begin{align} \label{eq: VI of tilde V}
    \begin{cases}
        \min\{-\widetilde{V}_t - \frac{\theta^2}{2}y^2 \widetilde{V}_{yy} - (\beta-r) y \widetilde{V}_y + \beta \widetilde{V}, \, \widetilde{V} - \widetilde{U}_K\} = 0, & \text{for } (t,y) \in Q_y, \\
        \widetilde{V}(T,y) = \widetilde{U}_K(y), & \text{for } y \in (0,+\infty).
    \end{cases}
\end{align}

The primal and dual value functions are related by the verification theorem:
\begin{theorem}[Verification Theorem]
\label{thm: verification theorem}
Let $\widetilde{V}$ be the solution to \eqref{eq: VI of tilde V} such that $\widetilde{V} \in C^1(Q_y) \cap C(\overline{Q}_y)$, $\widetilde{V} \in C^{1,2,2}(Q_y \backslash \partial C_y)$ with locally bounded derivatives near $\partial C_y$, $|\widetilde{V}(t,y)| \leq C(y^q+1)$ for some constant $C>0$ and $q>0$, $\widetilde{V}(t,\cdot)$ is strictly convex, $- \widetilde{V}_y(t,y) \to +\infty$ as $y \to 0$ and $-\widetilde{V}_y(t,y) \to \hat{K} \leq K$ for some positive constant $\hat{K}$.
Then the primal value function $V(t,x)$ and the dual value function $\widetilde{V}(t,y)$ are related by
\begin{align} \label{eq: dual transform}
V(t,x) = \inf_{y>0} \left\{ \widetilde{V}(t,y) + xy \right\}, \quad (t,x) \in Q_x.
\end{align}
Moreover, for $s \geq t$, the optimal stopping time is $\tau^* = \inf\{s \geq t \ | (s,Y_s) \in S_y\} \, \wedge \, T,$ the optimal wealth is 
\begin{equation} \label{eq: optimal wealth}
    X_s^* = -\widetilde{V}_y(s, Y_s),
\end{equation}
and the optimal control is
\begin{align}
\pi_s^* &= \frac{\theta}{\sigma} Y_s \widetilde{V}_{yy}(s,Y_s) \mathbbm{1}_{\{t \leq s \leq \tau^*\}}, \label{eq: optimal control}
\end{align}
with $Y_t = y^*$ where $y^*$ is the unique solution to the equation $\widetilde{V}_y(t, Y_t) + x = 0$
with $X_t=x$ and $\mathbbm{1}_{\{\cdot\}}$ denotes the indicator function.
\begin{proof}
    The proof is standard. For the sake of completeness, we provide the proof in Appendix \ref{appendix: proof of verification theorem}.
\end{proof}
\end{theorem}

\textbf{Transformed dual value function.} To further simplify the linear variational inequality of $\widetilde{V}(t,y)$, we do a simple transformation of variables. Define 
\begin{align}
    \label{def: variable transform}
z = \log y, \ \tau = \frac{\theta^2}{2}(T-t), \ v(\tau,z) = \widetilde{V}(t,y).
\end{align}
The transformed function $v(\tau,z)$ satisfies a linear parabolic variational inequality:
\begin{align}
    \label{eq: variational inequality for v}
    \begin{cases}
        \min\{v_\tau - v_{zz} + \kappa v_z + \rho v, v-\widetilde{U}_K(e^z)\} = 0, & \text{for } (\tau,z) \in Q_z:=(0, \frac{\theta^2 T}{2}) \times \bR^1, \\
        v(0,z) = \widetilde{U}_K(e^z), & \text{for } z \in \bR^1,
    \end{cases}
\end{align}
where $\rho = \frac{2 \beta}{\theta^2}$ and $\kappa = \frac{2r-2 \beta}{\theta^2}+1$. Define $g(z):=\widetilde{U}_K(e^z)$. Define the linear parabolic operator of $v(\tau,z)$ as $\cG[v] := v_\tau - v_{zz} + \kappa v_z + \rho v.$

\textbf{Alignment with the conditions on $U$.} Two types of utility functions will be considered:
\begin{enumerate}
    \item Power utility: $U(x) = \frac{x^\gamma}{\gamma}$ for $0<\gamma<1$. The dual function is $\widetilde{U}_K(y) = \frac{1-\gamma}{\gamma} y^{\frac{\gamma}{\gamma-1}} - Ky$.
    \item non-HARA utility: $U(x) = \frac{1}{3} H^{-3}(x) + H^{-1}(x) + xH(x)$, where $H(x) = \left(\frac{2}{-1+\sqrt{1+4x}}\right)^{\frac{1}{2}}$. The dual function is $\widetilde{U}_K(y) = \frac{1}{3} y^{-3} + y^{-1} - Ky$.
\end{enumerate}
Both utility functions satisfy the conditions above.

\section{Main methods} \label{sec:algorithms}
\subsection{Traditional algorithms}
There are several traditional algorithms to solve the mixed optimal stopping and control problem. As our method is based on the dual problem, we would compare with those which are also based on solving the dual problem first, and then use Theorem \ref{thm: verification theorem} to get the primal value and optimal control. Since the dual problem is purely an optimal stopping problem without control, the most common way is to use Binomial Tree Method (BTM), similarly to American option pricing. The other way is to find the Global Closed Approximation (GCA) of optimal stopping boundary, and then use Monte Carlo simulation to compute the dual value by stopping when the dual process touches it, see \citet{maGlobalClosedFormApproximation2019}. However, it is not always possible to find the closed form of boundary. Moreover, the number and shape of the stopping boundaries depend on the parameters of the model, and there is not a general representation for all cases, see Figure 1 in \citet{maGlobalClosedFormApproximation2019}. Both of the methods above require extensive simulations and are not efficient. Most importantly, they only give results for a single pair of $(t,x)$, and we need to repeat the same process even for a small change.

Another very important aspect of the problem is that to use the primal-dual relationship \eqref{eq: dual transform} and \eqref{eq: optimal control}, we need the derivatives of the dual value function, not only the values. \eqref{eq: dual transform} requires us to compute the first-order derivative $\widetilde{V}_y(t,y)$ and solve the equation $\widetilde{V}_y(t, Y_t) + x = 0$ to get $Y_t$. Based on it, \eqref{eq: optimal control} requires us to compute the second-order derivative $\widetilde{V}_{yy}(t,y)$ in extra. With the traditional methods above, we need to use FDM to compute the derivatives. As we said, the whole process needs to be repeated for any difference in $(t,x)$. Moreover, FDM is very sensitive to the grid size, and there is no benchmark for $V(t,x)$, $X_t^*$ and $\pi_t^*$, so it is difficult to choose the grid size in a convincing way. Furthermore, since we need to solve an equation involving the first-order derivative, we need to use a root-finding method, which requires FDM to be repeated many times and introduces extra errors.

\subsection{Deep neural network methods (DGM and FBR-DGM)} \label{sec:deep learning}

Instead of using the traditional methods, we can use deep learning to solve the mixed optimal stopping and control problem. The main idea is to use a neural network to approximate the solution to \eqref{eq: variational inequality for v}, and then transform it to $\widetilde{V}(t,y)$ as in \eqref{def: variable transform}. The advantage of it is that we get the dual value function on the whole domain $Q_y$. Hence, we can compute the derivatives $\widetilde{V}_y(t,y)$ and $\widetilde{V}_{yy}(t,y)$, as well as the primal value function $V(t,x)$ and the optimal control $\pi_t^*$ instantly once the training is done.

Since the domain $Q_z$ is unbounded, we need to truncate it to a bounded domain. We choose $z \in (-C_z, C_z)$ for some constant $C_z>0$. Define $\Omega=(-C_z,C_z)$, $\partial\Omega=\{-C_z,C_z\}$, $\Sigma=(0,T)\times\{-C_z,C_z\}$ and $\Omega_T=(0,T)\times \Omega$. The parabolic boundary is $\partial_p\Omega_T = \left(\{0\}\times [-C_z,C_z]\right) \cup \Sigma$, including the bottom $\{0\}\times\Omega$, the corner $\{0\}\times\partial\Omega$, and the lateral boundary $\Sigma$.

We use two different loss functions to train the neural network. The first one is a standard loss function, which is the mean square error (MSE) between the neural network output $f(t,z;\theta)$ and target functions: $J_1(f) = \|\min\{\cG[f], f-g\}\|^2_{L^2(\Omega_T), \mu_1} + \|f-g\|^2_{L^{2}(\partial_p\Omega_T), \mu_2},$
where $\mu_1, \mu_2$ are probability measures on $\Omega_T, \partial_p\Omega_T$ respectively. We require them to be absolutely continuous with respect to the Lebesgue measure and bounded away from zero. In practice, we can use uniform or truncated normal distribution. We use DGM to denote the network trained with the loss function $J_1$.

Based on the above, the second loss function includes the fractional boundary loss, which is the first time we propose in this paper. The idea of the fractional boundary loss is to match the boundary condition of the dual value function more than only $L^2$, which can improve the interior regularity of the solution. This is due to the inverse trace mapping from functional analysis. For instance, the most well-known result is the trace operator $\gamma: H^1(\Omega) \to H^{1/2}(\partial\Omega)$, where $H^1(\Omega)$ is the Sobolev space of functions with square integrable first-order derivatives and $H^{1/2}(\partial\Omega)$ is the Sobolev-Slobodeckij space of functions with square integrable fractional derivatives. $\gamma$ is surjective and has a right inverse,
$\cE : H^{1/2}(\partial\Omega) \to H^1(\Omega),$ such that $(\gamma \circ \cE)(g) = g$ for all $g \in H^{1/2}(\partial\Omega).$ Moreover, $\cE$ is continuous,
which ensures $\|\cE(x)\|_{H^1(\Omega)} \le C\,\|x\|_{H^{1/2}(\partial\Omega)}$ for a constant $C$. For the detailed definition of the Sobolev and Sobolev-Slobodeckij space, one can consult \citet{leoniFirstCourseSobolev2017} and \citet{brezisFunctionalAnalysisSobolev2011}. Intuitively speaking, there is a gain in interior regularity from boundary regularity. By incorporating the boundary norm in the loss function, we obtain a better convergence in the interior.

In our specific case, we need to consider the parabolic operator $\cG$. According to Theorem 2.3 in Chapter 4 of \citet{lionsNonHomogeneousBoundaryValue1972a}, we choose the following norms of the fractional Sobolev space on $\Sigma$, which are simplified since $\Omega$ is one-dimensional and tangential derivatives disappear:
\begin{align*}
    & \|d\|_{H^{\frac{3}{4}, \frac{3}{2}}(\Sigma)}\\
    =& \left(\|d\|_{L^2(\Sigma)}^2 + \int_0^T\int_0^T \frac{|d(t,C_z)-d(\tau,C_z)|^2+|d(t,-C_z)-d(\tau,-C_z)|^2}{|t-\tau|^{2.5}} \, dt \, d\tau\right)^{1/2},
\end{align*}
and
\begin{align*}
    \|d\|_{H^{\frac{1}{4}, \frac{1}{2}}(\Sigma)} =& \left(\|d\|_{L^2(\Sigma)}^2+\int_0^T \frac{|d(t,C_z)-d(t,-C_z)|^2}{4C_z^2} \, dt\right.\\
    &+ \left.\int_0^T\int_0^T \frac{|d(t,C_z)-d(\tau,C_z)|^2+|d(t,-C_z)-d(\tau,-C_z)|^2}{|t-\tau|^{1.5}} \, dt \, d\tau \right)^{1/2}.
\end{align*}
The fractional boundary loss function is defined as
\[
    J_2(f) = J_1(f) + \|f(0,\cdot)-g(0,\cdot)\|_{H^1(\Omega),\mu_3}^2+\|f-g\|^2_{H^{\frac{3}{4},\frac{3}{2}}(\Sigma),\mu_4} + \|f_z-g_z\|^2_{H^{\frac{1}{4},\frac{1}{2}}(\Sigma),\mu_4},
\]
where $\mu_3$ and $\mu_4$ are the probability measures on $\Omega$ and $\Sigma$. $f_z$ and $g_z$ are the derivatives of $f$ and $g$ with respect to $z$. We use FBR-DGM to denote the network trained with the loss function $J_2$.

\section{Numerical experiments} \label{sec:numerical experiments}

\subsection{Comparison among approximations of \texorpdfstring{$V$}{optimal value}} \label{sec:numerical 1}

Following Section \ref{sec:algorithms}, we check the neural network performance for primal value function by comparing the approximations with different methods. For the detailed algorithm and implementation, see Appendix \ref{appendix:algorithm 1}. 
When applying Theorem \ref{thm: verification theorem}, we need to compute the minimizer of \eqref{eq: dual transform} for each $(t,x)$. We can use bisection method or grid search. 
Bisection means we choose the minimizer for \eqref{eq: dual transform} from bisection method towards $\partial_y(\widetilde{V}(t,y) + x y)$. Grid means we evaluate the value for \eqref{eq: dual transform} on a dense grid of 200 points of $y$ and choose the minimizer. 
Numerical experiments show that the two methods are very similar, reflecting that the neural network converges up to first-order derivative. We only report the results with Grid.
The first method to compare is BTM. The other one is GCA from \citet{maGlobalClosedFormApproximation2019}, which approximates the free boundary for $\widetilde{V}$, and then uses FDM and bisection method to solve the first-order condition.

\begin{example}
\label{example: 1}
For power and non-HARA utility, we compare the numerical results among DGM, FBR-DGM, GCA and BTM. The parameters used are $\mu = 0.1$, $\beta = 0.1$, $r = 0.05$, $\sigma = 0.3$, $K = 1$, $\gamma = 0.5$, $T = 1$, $t = 0$, initial wealth $x_0 \in [1.1, 1.2, 1.3, 1.4, 1.5, 1.6, 1.7, 1.8, 1.9, 2.0]$, number of time steps for BTM $N = 2000$.
\end{example}

\begin{figure}[htbp]
  \centering
  \includegraphics[width=\textwidth]
                   {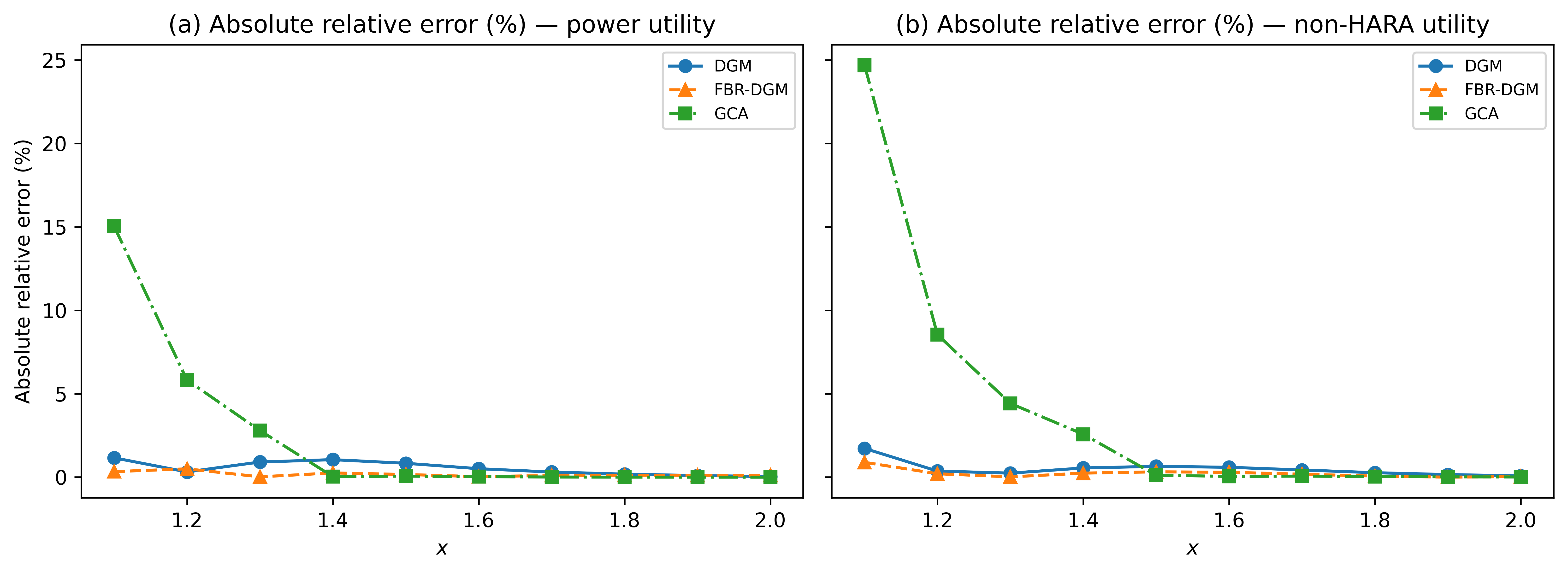}
  \caption{Approximation quality of DGM, FBR-DGM, and GCA
           compared with BTM.  Panels (a)-(b) show the 
           absolute relative errors for two utilities,
           averaged over five training seeds.}
  \label{figure:1}
\end{figure}

Figure \ref{figure:1} and Table \ref{tab1} illustrates the comparisons in Example \ref{example: 1}, benchmarked with BTM.

\begin{table}[htbp]
    \centering
    \begin{threeparttable}
    \caption{Comparison of methods, benchmarked with BTM, averaged over seeds.}
    \label{tab1}
    \begin{tabular}{llcccc}
        \toprule
        Utility & Method
        & \makecell{Mean Abs.\\Rel. Diff. (\%)\tnote{\S}}
        & \makecell{Std. Rel.\\Diff. (\%)\tnote{\S}}
        & \makecell{Offline\\ Train (s)}
        & \makecell{Online\\ /Pair (ms)\tnote{\S}}\\
        \midrule
        \multirow{3}{*}{Power}
        & DGM & $0.547 \pm 0.077$ \tnote{\P} & $0.648\pm0.026$ & \best{599\pm11} & \best{22\pm0} \\
        & FBR-DGM\tnote{\dag} & \best{0.263\pm0.077} & \best{0.407\pm0.169} & $846\pm5$ & $21\pm0$ \\
        & GCA & $2.379$ & $4.856$ & - & $19931$ \\
        \midrule
        \multirow{3}{*}{non-HARA} 
        & DGM & $0.628\pm0.126$ & $0.753\pm0.297$ & \best{577\pm3} & $18\pm0$ \\
        & FBR-DGM\tnote{\dag} & \best{0.314\pm0.127} & \best{0.420\pm0.122} & $847\pm8$ & $18\pm0$ \\
        & GCA & 4.051 & 7.797 & - & 32551 \\
        \bottomrule
    \end{tabular}
    \begin{tablenotes}
        \item[\dag] Our proposed method.
        \item[\S] Computed across 10 spatial points, then 5 independent seeds for training the neural network.
        \item[\P] Entries are $\mu_{seed} \pm \sigma_{seed}$, i.e., mean $\pm$ std (1-$\sigma$) across 5 independent seeds for training.
    \end{tablenotes}
\end{threeparttable}
\end{table}

As we can see from above, the neural network methods are much more accurate than GCA. There is a noticeable improvement when including the fractional boundary loss, for both utilities in terms of mean and standard deviation, and the training time does not increase significantly.
The offline training time of neural networks is longer than the total evaluation time of GCA at 10 pairs of points, but it is still acceptable. Most importantly, if we already trained the network, it is instant to get $V(t,x)$ for different pairs of $(t,x)$. Instead, the whole processes of both GCA and BTM need to be repeated for different pairs of $(t,x)$, so the computation time will grow linearly. This is a huge advantage of the neural network method, and allows us to test it further with methods in the next sections.
\subsection{Consistency between \texorpdfstring{$V$}{optimal value} and \texorpdfstring{$X^*$}{optimal wealth}} \label{sec:numerical 2}
In Section \ref{sec:numerical 1}, we use the primal-dual relationship, \eqref{eq: dual transform}, directly to get the primal value $V(0,x_0)$ for given pairs. Based on Theorem \ref{thm: verification theorem}, $V(0,x_0)$ can also be obtained from the relation \eqref{eq: optimal wealth}.
We would compare it with the primal value that we got in Section \ref{sec:numerical 1}, and get the mean and standard deviation of the relative difference, which can verify if the neural network method admits the verification theorem and is self-consistent, see Algorithm \ref{alg:self_consistency_checks}. For the details, see Algorithm \ref{alg:self_consistency_optimal_wealth} in Appendix \ref{appendix:self_consistency_checks}.
\begin{algorithm}[htbp]
\caption{Self-consistency checks in Section \ref{sec:numerical 2} and Section \ref{sec:numerical 3}}
\label{alg:self_consistency_checks}
\begin{algorithmic}[1]
\For{$i = 1,2,\dots,M$}
  \State Simulate $\{Y_{t_n}^{(i)}\}_{n=0}^N$ via \eqref{eq: dual process}.
  \State Compute $\tau^{*(i)}
        = \min\{t_n:\,\widetilde U_K(Y_{t_n}^{(i)}) \ge \widetilde V(t_n,Y_{t_n}^{(i)})\}\;\wedge\;T$.
  \State {\bfseries(a) [Sec.~\ref{sec:numerical 2}]} Evaluate the stopped primal wealth $X_{\tau^{*(i)}}^{(i)} = -\,\widetilde V_{y}\bigl(\tau^{*(i)},Y_{\tau^{*(i)}}^{(i)}\bigr)$ via \eqref{eq: optimal wealth}.
  \State {\bfseries(b) [Sec.~\ref{sec:numerical 3}]} Simulate $X_{\tau^{*(i)}}^{\pi^*,(i)}$ via \eqref{eq: optimal control} with $\pi^*_t \;=\;\frac{\theta}{\sigma}\;Y_t\,\widetilde V_{yy}(t,Y_t)$.
  \State Compute discounted utility $P_i = e^{-\beta\,\tau^{*(i)}}\,U\bigl(X_{\tau^{*(i)}}^{(i)} - K\bigr)$.
  \State Compute relative error $r_i \;=\;\frac{P_i - V_0}{V_0}$, where $V_0=V(0,x_0)$ from Section \ref{sec:numerical 1}.
\EndFor
\State Compute sample mean $\bar r$ and sample standard deviation $s_r$ as
\begin{align*}
    \bar r = \frac{1}{M}\sum_{i=1}^M r_i, \quad s_r = \sqrt{\frac{1}{M-1}\sum_{i=1}^M (r_i - \bar r)^2}.
\end{align*}
\end{algorithmic}
\end{algorithm}

\begin{figure}[htbp]
    \centering
    \begin{subfigure}[b]{0.45\linewidth}
        \centering
        \includegraphics[width=\linewidth]{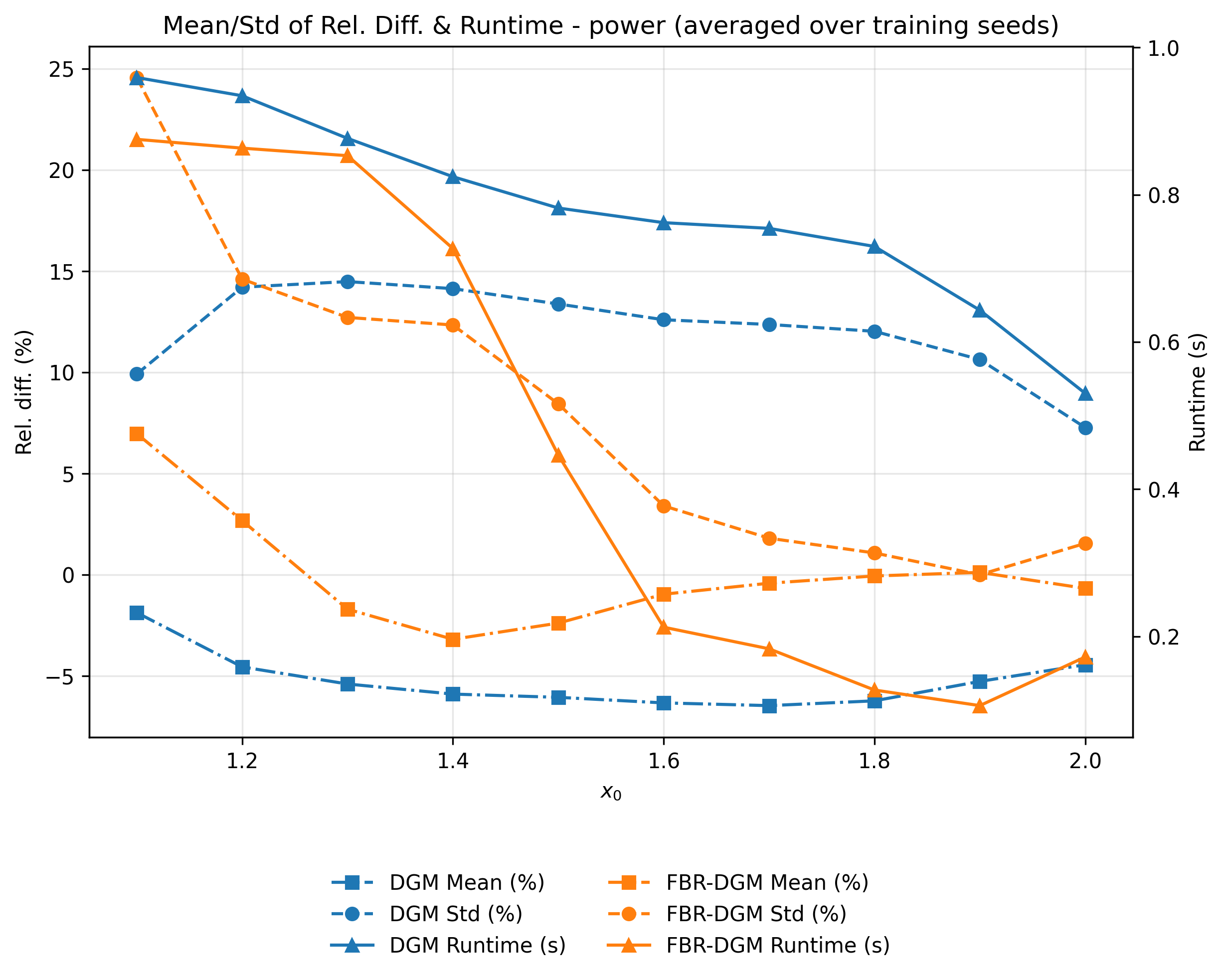}
        \caption{Power utility}
    \end{subfigure}
    \hfill
    \begin{subfigure}[b]{0.45\linewidth}
        \centering
        \includegraphics[width=\linewidth]{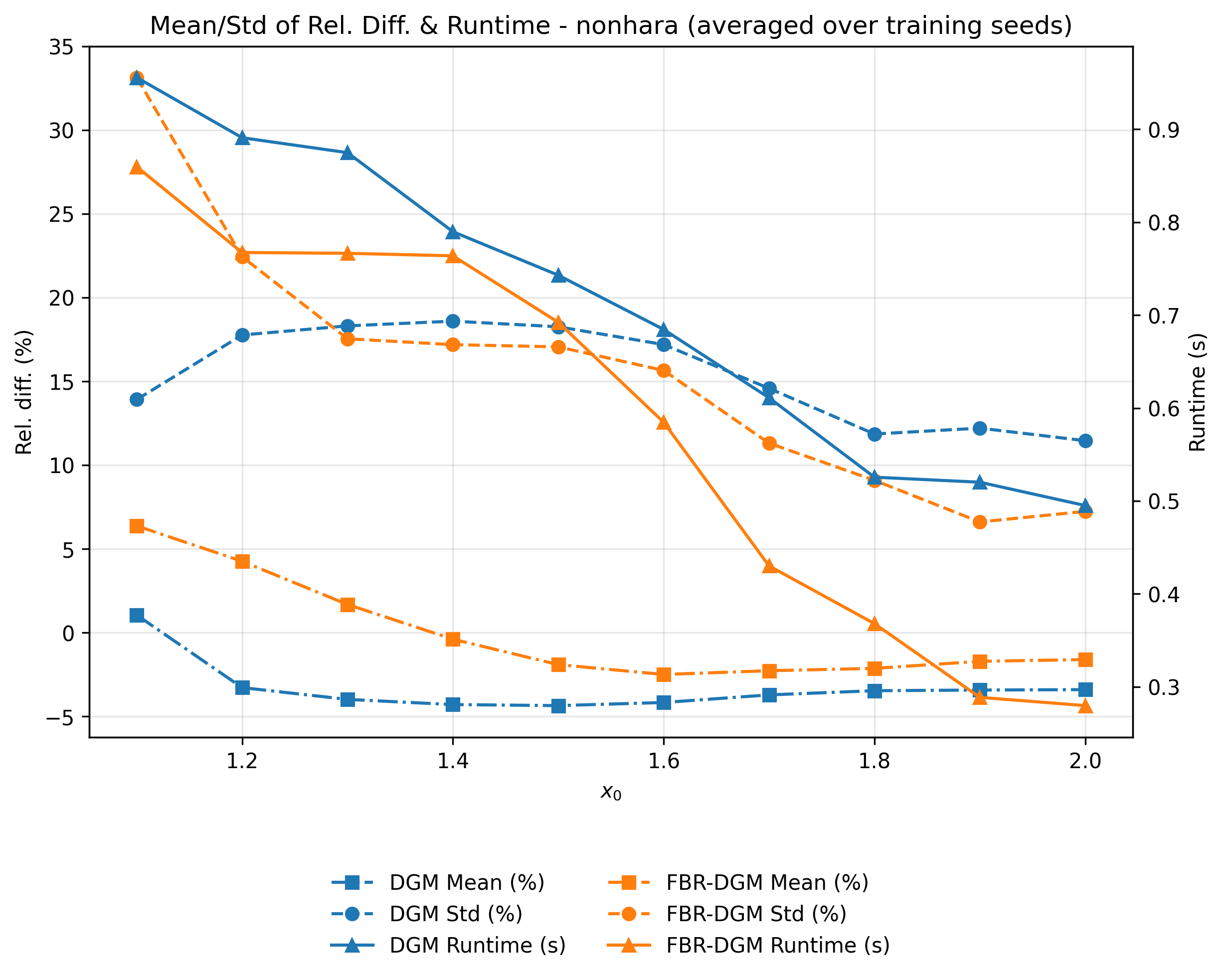}
        \caption{non-HARA utility}
    \end{subfigure}

    \caption{Results of Example \ref{example: 2} for networks trained with DGM and FBR-DGM, averaged over seeds.}
    \label{fig:2}
\end{figure}

\begin{example}
    \label{example: 2}
    With the same parameters as Example \ref{example: 1}, we compare the relative difference between $e^{-\beta \tau^*}U(X_{\tau^*} - K)$ from \eqref{eq: optimal wealth} and $V(0, x_0)$ for power utility and non-HARA utility with $M$ independent simulations. We choose $\Delta t = 0.01$ and $M = 100$.
\end{example}
Figure \ref{fig:2} and Table \ref{tab2} illustrate the results in Example \ref{example: 2}. For reliability, we use 5 different random seeds for the Monte-Carlo simulation, and the results are averaged over the seeds at each $x_0$. As we can see, the mean of relative difference is typically within $5\%$, and it decreases as $x_0$ increases, which is consistent with Figure \ref{figure:1}, where we can infer that the continuation region is smaller for greater $x_0$ by comparing the immediate exercise value and the primal value.

\begin{table}[htbp]
    \centering
    \begin{threeparttable}
    \caption{Results of Example \ref{example: 2} for networks trained with DGM and FBR-DGM}
    \label{tab2}
    \begin{tabular}{llccc}
        \toprule
        Utility & Method & Mean Abs. Rel. Diff. (\%)\tnote{\S} & Std Rel. Diff. (\%)\tnote{\S} & Average Time (s)\tnote{\S} \\
        \midrule
        \multirow{2}{*}{Power}
        & DGM & $5.258\pm0.807$\tnote{\P} & $12.100\pm0.664$ & $0.780\pm0.093$ \\
        & FBR-DGM\tnote{\dag} & \best{2.006\pm0.427} & \best{8.049\pm1.126} & \best{0.457\pm0.040} \\
        \midrule
        \multirow{2}{*}{non-HARA} 
        & DGM & $3.637\pm1.857$ & \best{15.412\pm4.502} & $0.709\pm0.197$ \\
        & FBR-DGM\tnote{\dag} & \best{2.776\pm1.367} & $15.723\pm3.155$ & \best{0.580\pm0.118} \\
        \bottomrule
    \end{tabular}
    \begin{tablenotes}
        \item[\dag] Our proposed method.
        \item[\S] Computed across 5 seeds for simulation, 10 spatial points, then 5 seeds for training.
        \item[\P] Entries are $\mu_{seed} \pm \sigma_{seed}$, i.e., mean $\pm$ std (1-$\sigma$) across 5 independent seeds for training.
    \end{tablenotes}
\end{threeparttable}
\end{table}

\subsection{Consistency between \texorpdfstring{$V$}{optimal value} and \texorpdfstring{$\pi^*$}{optimal control}} \label{sec:numerical 3}
Apart from the approach in Section \ref{sec:numerical 2}, there is another way to verify the self-consistency of the neural network from the optimal control in \eqref{eq: optimal control}, see Algorithm \ref{alg:self_consistency_checks}. 
For the details, see Appendix \ref{appendix:self_consistency_checks}.

\begin{example}
    \label{example: 3}
With the same parameters as Example \ref{example: 1}, we compare the relative difference between $e^{-\beta \tau^*}U(X_{\tau^*} - K)$ from \eqref{eq: optimal control} and $V(0, x_0)$ for power utility and non-HARA utility with $M$ independent simulations. We choose $\Delta t = 0.01$ and $M = 100$.
\end{example}

\begin{figure}[htbp]
    \centering
    \begin{subfigure}[b]{0.45\linewidth}
        \centering
        \includegraphics[width=\linewidth]{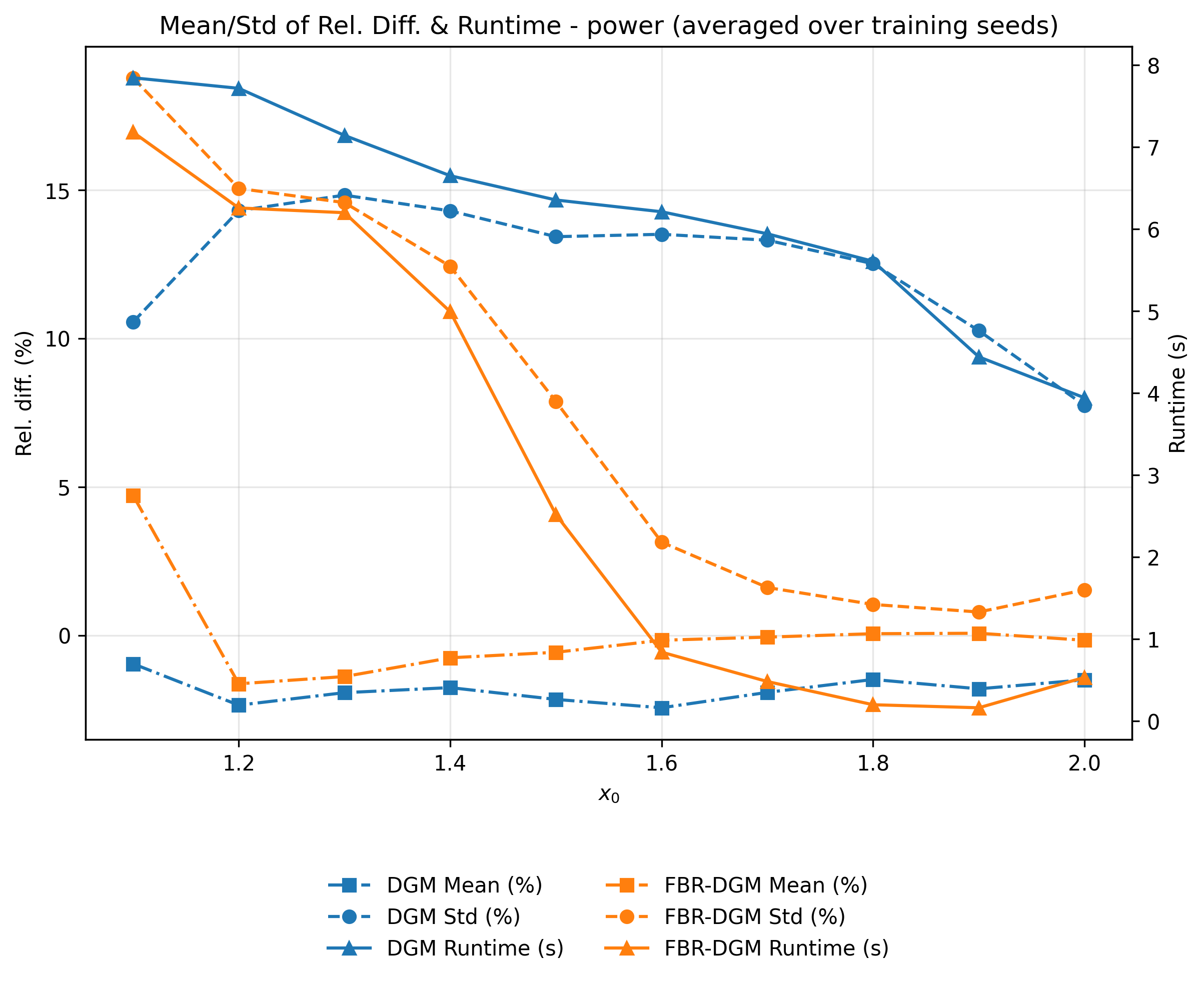}
        \caption{Power utility}
    \end{subfigure}
    \hfill
    \begin{subfigure}[b]{0.45\linewidth}
        \centering
        \includegraphics[width=\linewidth]{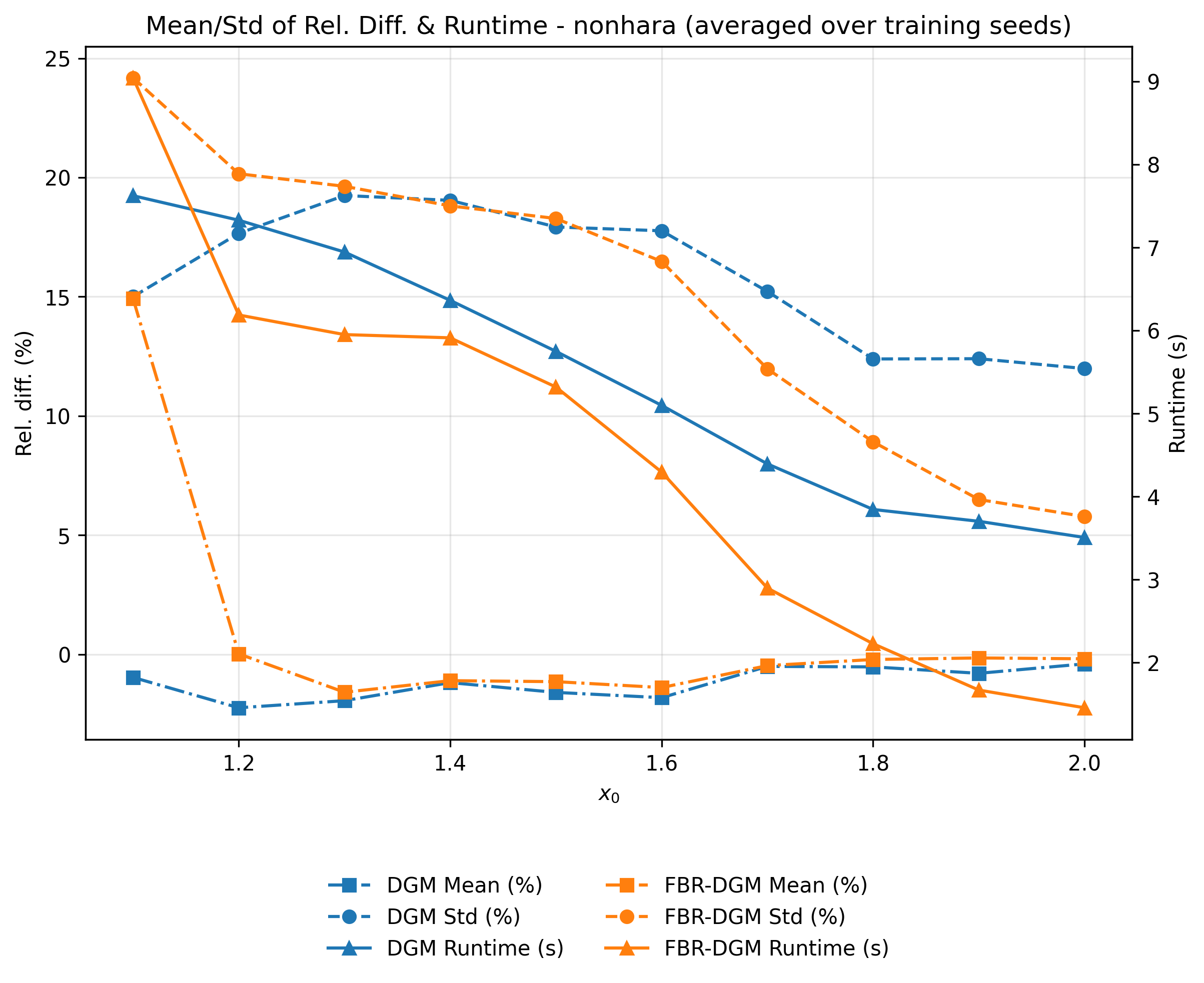}
        \caption{non-HARA utility}
    \end{subfigure}

    \caption{Results of Example \ref{example: 3} for networks trained with DGM and FBR-DGM, averaged over seeds.}
    \label{fig:3}
\end{figure}

Figure \ref{fig:3} and Table \ref{tab3} illustrate the results in Example \ref{example: 3}. As we can see, the performance of FBR-DGM is better for power utility than for non-HARA utility, and it is more consistent with a higher $x_0$.

\begin{table}[htbp]
    \centering
    \begin{threeparttable}
    \caption{Results of Example \ref{example: 3} for networks trained with DGM and FBR-DGM}
    \label{tab3}
    \begin{tabular}{llccc}
        \toprule
        Utility & Method & Mean Abs. Rel. Diff. (\%)\tnote{\S} & Std Rel. Diff. (\%)\tnote{\S} & Average Time (s)\tnote{\S} \\
        \midrule
        \multirow{2}{*}{Power}
        & DGM & $1.832\pm0.306$\tnote{\P} & $12.478\pm0.938$ & $6.187\pm0.885$ \\
        & FBR-DGM\tnote{\dag} & \best{1.044\pm0.422} & \best{7.684\pm0.914} & \best{2.937\pm0.309} \\
        \midrule
        \multirow{2}{*}{non-HARA} 
        & DGM & \best{1.247\pm0.443} & $15.862\pm4.865$ & $5.455\pm1.860$ \\
        & FBR-DGM\tnote{\dag} & $2.315\pm0.681$ & \best{15.069\pm3.462} & \best{4.495\pm1.093} \\
        \bottomrule
    \end{tabular}
    \begin{tablenotes}
        \item[\dag] Our proposed method.
        \item[\S] Computed across 5 seeds for simulation, 10 spatial points, then 5 seeds for training.
        \item[\P] Entries are $\mu_{seed} \pm \sigma_{seed}$, i.e., mean $\pm$ std (1-$\sigma$) across 5 independent seeds for training.
    \end{tablenotes}
\end{threeparttable}
\end{table}

\section{Conclusion} \label{sec:conclusion}
In this paper, we propose a deep learning method to solve the mixed optimal stopping and control problem. We use a neural network to approximate the dual value function, which is the solution to a variational inequality, and then transform it to the primal value function. To improve the convergence and accuracy of the derivatives, we propose a new loss function that includes the fractional boundary norms, theoretically justified by the inverse trace mapping. Taking advantage of the primal-dual relationship and the efficiency of the neural network, we design two innovative algorithms to verify the self-consistency of the neural network. The numerical results show that our method is more accurate and efficient than traditional methods. Compared with the classical $L^2$ loss function, the new loss function can significantly improve the accuracy and stability of the results.

There are several directions for future work. First, we can extend our method to more general settings, such as partial information, multiple assets, and transaction costs. Second, we can explore more theoretical properties of the Fractional-Boundary-Regularized loss function, such as convergence rates and regularity.

\bibliography{references}

\newpage
\appendix
\section{Proof of Theorem \ref{thm: verification theorem}} \label{appendix: proof of verification theorem}
\begin{proof}
    We only show the proof for $t=0$ and $X_0=x \in (K,+\infty)$. For a general $t \in [0,T]$ and $X_t \in (K,+\infty)$, the proof is similar.

    Under the conditions on $\widetilde{V}$, the following function is well-defined:
    \[
    \varphi(t,x) := \inf_{y>0} \{\widetilde{V}(t,y)+xy\}, \quad (t,x) \in Q_x.
    \]
    Since $\widetilde{V} \in C^1(Q_y) \cap C(\overline{Q}_y)$ and $\widetilde{V}(t,\cdot)$ is strictly convex, from duality,
    \[
    \widetilde{V}(t,y) = \sup_{x>K} \{\varphi(t,x)-xy\}, \quad (t,y) \in Q_y.
    \]
    Since $- \widetilde{V}_y(t,y) \to +\infty \text{ as } y \to 0 \text{ and } -\widetilde{V}_y(t,y) \to \hat{K} \leq K$ for some positive constant $\hat{K}$, by first-order condition, there is a unique $y^\star(t,x) \in (0,+\infty)$ such that
    \[
    \widetilde{V}_y(t,y^\star(t,x)) + x = 0, \quad \varphi(t,x) = \widetilde{V}(t,y^\star(t,x)) + xy^\star(t,x).
    \]
    Differentiate the above equation of $\varphi(t,x)$ with respect to $x$,
    \[
    \varphi_x(t,x) = \widetilde{V}_y(t,y^\star(t,x)) y_x^\star(t,x) + x y_x^\star(t,x) + y^\star(t,x) = y^\star(t,x).
    \]
    Hence,
    \[
    \widetilde{V}(t,\varphi_x(t,x)) = \varphi(t,x)-x\varphi_x(t,x), \quad \widetilde{V}_y(t,\varphi_x(t,x)) = -x.
    \]
    Differentiate the above equations of $\widetilde{V}$ with respect to $t$ and $\widetilde{V}_y$ with respect to $x$,
    \[
    \widetilde{V}_t(t,\varphi_x(t,x)) = \varphi_t(t,x), \quad \widetilde{V}_{yy}(t,\varphi_x(t,x)) = -\frac{1}{\varphi_{xx}(t,x)}.
    \]
    From the above equations,
    \begin{align} \label{eq: ieq1}
    & -\varphi_t + \frac{\theta^2}{2} \frac{\varphi_x^2}{\varphi_{xx}} - r\varphi_x x + \beta \varphi \\
    = & -\widetilde{V}_t - \frac{\theta^2}{2}y^{*2} \widetilde{V}_{yy} - (\beta-r) y^* \widetilde{V}_y + \beta \widetilde{V} \geq 0 \nonumber,
    \end{align}
    where the inequality follows from \eqref{eq: VI of tilde V}. Similarly, since $U(x-K) = \inf_{y>0} \{\widetilde{U}_K(y)+xy\}$ by applying the duality on \eqref{def: dual utility function},
    \begin{equation} \label{eq: ieq2}
    \varphi(t,x)-U(x-K) = \inf_{y>0} \{\widetilde{V}(t,y)+xy\} - \inf_{y>0} \{\widetilde{U}_K(y)+xy\} \geq 0,
    \end{equation}
    where the inequality follows from $\widetilde{V}(t,y)\geq \widetilde{U}_K(y)$ in \eqref{eq: VI of tilde V}. From the same argument, there are only two cases in $Q_x$: $\widetilde{V}(t,y^*(t,x)) > \widetilde{U}_K(y^*(t,x))$ and $\widetilde{V}(t,y^*(t,x)) = \widetilde{U}_K(y^*(t,x))$. When $\widetilde{V}(t,y^*(t,x)) > \widetilde{U}_K(y^*(t,x))$, the other term in \eqref{eq: VI of tilde V} is equal to zero, and
    \begin{align} \label{eq: eq1}
    & -\varphi_t + \frac{\theta^2}{2} \frac{\varphi_x^2}{\varphi_{xx}} - r\varphi_x x + \beta \varphi \\
    = & -\widetilde{V}_t - \frac{\theta^2}{2}y^{*2} \widetilde{V}_{yy} - (\beta-r) y^* \widetilde{V}_y + \beta \widetilde{V} = 0 \nonumber.
    \end{align}
    When $\widetilde{V}(t,y^*(t,x)) = \widetilde{U}_K(y^*(t,x))$, the function $\widetilde{V}(t,y)-\widetilde{U}_K(y)$ attains its minimum at $y^*(t,x)$, and the first-order condition is satisfied. Hence,
    \[
    \widetilde{V}_y(t,y^*(t,x)) = \widetilde{U}_K'(y^*(t,x)) = -x.
    \]
    Hence, $y^*$ is a stationary point of $\widetilde{U}_K(y)+xy$. Since $\widetilde{U}_K$ is strictly convex,
    \[
    \widetilde{U}_K(y^*(t,x)) + xy^*(t,x) = \inf_{y>0} \{\widetilde{U}_K(y) + xy\} = U(x-K).
    \]
    Hence, in this case,
    \begin{equation} \label{eq: eq2}
    U(x-K)=\widetilde{V}(t,y^*(t,x))+x y^*(t,x) = \varphi(t,x).
    \end{equation}
    Combining the inequalities \eqref{eq: ieq1} and \eqref{eq: ieq2}, and the equalities \eqref{eq: eq1} and \eqref{eq: eq2}, we have
    \begin{equation} \label{eq: VI of varphi}
    \min\{-\varphi_t + \frac{\theta^2}{2} \frac{\varphi_x^2}{\varphi_{xx}} - r\varphi_x x + \beta\varphi, \varphi(t,x)-U(x-K)\} = 0, \quad (t,x) \in Q_x.
    \end{equation}
    Comparing \eqref{eq: VI of varphi} with \eqref{eq: VI of V without control}, we can see $\varphi$ satisfies the same variational inequality as $V$. Equivalently, $\varphi$ satisfies \eqref{eq: VI of V with control}. Hence,
    \begin{align*}
        & d(e^{-\beta t} \varphi(t,X_t)) \\
        = & e^{-\beta t}\left(-\beta \varphi(t,X_t)+\varphi_t(t,X_t)+\varphi_x(t,X_t)(rX_t+\sigma\theta\pi_t)+\frac{\sigma^2}{2}\varphi_{xx}(t,X_t)\pi_t^2\right)dt \\
        & + \sigma e^{-\beta t}\varphi_{x}(t,X_t)\pi_t dW_t \\
        \leq & \sigma e^{-\beta t}\varphi_{x}(t,X_t)\pi_t dW_t,
    \end{align*}
    where the inequality follows from \eqref{eq: VI of V with control}. Hence, $(e^{-\beta t}\varphi(t,X_t))_{0\leq t\leq T}$ is a supermartingale. For any stopping time $\tau \in [0,T]$, define
    \[
    \tau_n := \inf\left\{0\leq t\leq T \, | \, X_t \geq x+n \text{ or } |\pi_t|\geq n\right\} \wedge \tau.
    \]
    Then,
    \[
    E[e^{-\beta \tau_n} \varphi(\tau_n,X_{\tau_n})] \leq \varphi(0,x).
    \]
    From Fatou's lemma,
    \[
    \varphi(0,x) \geq \lim_{n\to\infty} E[e^{-\beta\tau_n} \varphi(\tau_n,X_{\tau_n})] \geq E[e^{-\beta\tau}\varphi(\tau,X_\tau)] \geq E[e^{-\beta\tau}U(X_\tau-K)].
    \]
    Taking arbitrary $\pi,\tau$ in the above inequality,
    \[
    \varphi(0,x) \geq \sup_{\pi,\tau} E[e^{-\beta\tau}U(X_\tau-K)] = V(0,x).
    \]
    Define
    \[
    \pi^* = -\frac{\theta \varphi_x}{\sigma \varphi_{xx}}, \quad \tau^* = \inf\{t\geq 0 \, | \, \varphi(t,X_t)=U(X_t-K)\} \wedge T.
    \]
    Then,
    \begin{align*}
        & d(e^{-\beta t} \varphi(t,X_t^{\pi^*})) \\
        = & e^{-\beta t}\left(-\beta \varphi(t,X_t^{\pi^*})+\varphi_t(t,X_t^{\pi^*})+\varphi_x(t,X_t^{\pi^*})(rX_t^{\pi^*}+\sigma\theta\pi_t^*)+\frac{\sigma^2}{2}\varphi_{xx}(t,X_t^{\pi^*})\pi_t^{*2}\right)dt \\
        & + \sigma e^{-\beta t}\varphi_{x}(t,X_t^{\pi^*})\pi_t^* dW_t \\
        = & e^{-\beta t} \left(-\beta\varphi(t,X_t^{\pi^*})+\varphi_t(t,X_t^{\pi^*})+rX_t^{\pi^*}\varphi_x(t,X_t^{\pi^*})-\frac{\theta^2\varphi_x^2(t,X_t^{\pi^*})}{2\varphi_{xx}(t,X_t^{\pi^*})}\right)dt \\
        & +\sigma e^{-\beta t}\varphi_{x}(t,X_t^{\pi^*})\pi_t^* dW_t \\
        = & \sigma e^{-\beta t}\varphi_{x}(t,X_t^{\pi^*})\pi_t^* dW_t, \quad t\in [0,\tau^*],
    \end{align*}
    where the last equality follows from \eqref{eq: VI of varphi}. Integrating the above equation from $0$ to $\tau^*$ and taking expectations, we have
    \[
    \varphi(0,x) = E[e^{-\beta \tau^*} \varphi(\tau^*,X_{\tau^*}^{\pi^*})] = E[e^{-\beta\tau^*}U(X_{\tau^*}^{\pi^*}-K)] \leq V(0,x).
    \]
    Hence,
    \[
    \varphi(0,x)=V(0,x).
    \]
    We also conclude that $\tau^*$ and $\pi^*$ are optimal stopping time and control, respectively. From previous arguments, the optimal wealth $x=-\widetilde{V}_y(t,y^*(t,x))$.
\end{proof}

\section{Algorithms} \label{appendix:algorithms}
\subsection{Algorithm for Section \ref{sec:numerical 1} and implementation details} \label{appendix:algorithm 1}
The algorithm for the primal value function in Section \ref{sec:numerical 1} is given in Algorithm \ref{alg:primal_from_network}. The implementation details for the neural network are given in Implementation \ref{impl:nn_details}.

\begin{algorithm}[htbp]
\caption{Primal Value from Neural Network Approximation: 
  $v_n(\tau,z) \rightarrow \widetilde{V}_n(t,y) \rightarrow V_n(t,x)$}
\label{alg:primal_from_network}
\begin{algorithmic}[1]
\Require Model parameters: $\theta, T, K, \kappa, \rho$; constants for domain: $C_z, C_x^1, C_x^2 > 0$; utility function $U$.
\Ensure Approximate primal value function $V_n(t,x)$ for $(t,x) \in (0,T) \times (C_x^1, C_x^2)$.

\State \textbf{Step 1: Train Neural Network $v_n(\tau,z)$}
\begin{itemize}
    \item \textbf{Domain:}
    \[
        Q := (-C_z, C_z), \qquad Q_T := \Bigl(0,\tfrac{\theta^2T}{2}\Bigr)\times Q.
    \]
    \item \textbf{Sampling:} Draw training points $(\tau,z)$ from a truncated normal distribution centered at
    \[
        \left(\tfrac{\theta^2T}{4},\,0\right).
    \]
    \item \textbf{Training:} Learn a neural network approximation $v_n(\tau,z)$ to the solution of the variational inequality
    \[
        \min\{v_\tau - v_{zz} + \kappa\,v_z + \rho\,v,\, v - g\} = 0 \quad \text{in } Q_T,
    \]
    subject to boundary conditions
    \[
        v = g \quad \text{on } \partial_p Q_T,
    \]
    where
    \[
        g(z) := \widetilde{U}_K(e^z), \quad 
        \partial_p Q_T := \left(\{0\} \times [-C_z, C_z]\right)
        \cup \left(\left(0,\tfrac{\theta^2 T}{2}\right) \times \{-C_z, C_z\}\right).
    \]
\end{itemize}

\State \textbf{Step 2: Transform to Dual Value $\widetilde{V}_n(t,y)$}
\begin{itemize}
    \item \textbf{Change of Variables:}
    \[
        t = T - \frac{2\tau}{\theta^2}, \qquad y = e^z.
    \]
    \item \textbf{Output:} Define the dual value function by
    \[
        \widetilde{V}_n(t,y) := v_n\left(\tfrac{\theta^2}{2}(T - t), \ln y\right),
    \]
    defined on $(0,T) \times (e^{-C_z}, e^{C_z})$.
\end{itemize}

\State \textbf{Step 3: Recover Primal Value $V_n(t,x)$}
\begin{itemize}
    \item \textbf{Dual-Primal Transformation:}
    \[
        V_n(t,x) := \min_{y \in (e^{-C_z}, e^{C_z})} \left\{ \widetilde{V}_n(t,y) + x y \right\}.
    \]
    \item \textbf{Output:} $V_n(t,x)$ is defined on $(0,T) \times (C_x^1, C_x^2)$.
\end{itemize}

\end{algorithmic}
\end{algorithm}

\begin{implementation}[htbp]
\caption{Neural Network Implementation Details}
\label{impl:nn_details}
\begin{algorithmic}[1]
\State \textbf{Activation Function:} Use the smooth \textbf{Mish} activation,
\[
  \mathrm{Mish}(x) \;=\; x \,\tanh\left(\ln(1+e^x)\right),
\]
to guarantee $C^2$-smoothness for the parabolic operator and universal approximation theorem.

\State \textbf{Optimizer:} Use the \textbf{Adam} optimizer with default momentum parameters.

\State \textbf{Domain Scaling:}
\begin{itemize}
  \item \emph{Training domain:} $z\in(-C_z,C_z)$ with $C_z=1.5$.
  \item \emph{Evaluation domain:} $z\in(-C_z^{\mathrm{eval}},C_z^{\mathrm{eval}})$ with $C_z^{\mathrm{eval}}=1.0$ to ensure stable derivative estimates.
\end{itemize}

\State \textbf{Loss Function:} Based on the definitions in Section \ref{sec:deep learning}, we would use $J_1$ for DGM and $J_2$ for FBR-DGM. We scale the loss function by $1+|g(z)|$ or $1+|g'(z)|$ to counteract large variations in $g(z)$. The explicit loss functions are given by
\[
J_1^{\mathrm{scale}}(f) = \left\|\frac{\min\{\cG[f], f-g\}}{1+|g|}\right\|^2_{L^2(\Omega_T), \mu_1} + \left\|\frac{f-g}{1+|g|}\right\|^2_{L^{2}(\partial_p\Omega_T), \mu_2}
\]
for DGM, and
\begin{align*}
J_2^{\mathrm{scale}}(f) &= J_1^{\mathrm{scale}}(f) \\ 
& + \left\|\frac{f(0,\cdot)-g(0,\cdot)}{1+|g(0,\cdot)|}\right\|_{H^1(\Omega),\mu_3}^2 + \left\|\frac{f-g}{1+|g|}\right\|^2_{H^{\frac{3}{4},\frac{3}{2}}(\Sigma),\mu_4} + \left\|\frac{f_z-g_z}{1+|g_z|}\right\|^2_{H^{\frac{1}{4},\frac{1}{2}}(\Sigma),\mu_4},
\end{align*}
for FBR-DGM.

\State \textbf{Dual Primal Minimization:} Solve
\[
  \min_{y\in(e^{-C_z},e^{C_z})}\bigl\{\widetilde V_n(t,y) + x\,y\bigr\}
\]
using either a bisection method or a direct grid search.

\State \textbf{Learning Rate Schedule:} Initialize at $2\times10^{-3}$ and halve every $4{,}000$ epochs over a total of $20{,}000$ epochs.

\State \textbf{Network Architecture:} Build an 8-layer feed-forward network with 128 neurons per layer, each followed by the Mish activation.

\State \textbf{Batch Sizes:}
\begin{itemize}
  \item Interior of $Q_T$: 5096 points,
  \item Initial-time boundary: 2048 points,
  \item Each lateral boundary side: 1024 points.
\end{itemize}

\State \textbf{Sampling Distribution:} Draw $\tau$ from a uniform distribution on $(0,\frac{\theta^2T}{2})$. For the loss from $J_1$, draw $z$ from a truncated normal centered at $0$ with standard deviation chosen so that $\mu + 3\sigma \approx C_z$. For the extra loss from $J_2$, draw $z$ from a uniform distribution on $(-C_z,C_z)$.

\State \textbf{Computational Environment:} All trainings were implemented in Python~3.11.12 with PyTorch~2.6.0 (CUDA~12.4), on a Google Colaboratory instance equipped with an NVIDIA A100 GPU (40 GB RAM) and 83.5 GB of system RAM.
\end{algorithmic}
\end{implementation}

\subsection{Full Algorithms for Section \ref{sec:numerical 2} and Section \ref{sec:numerical 3}} \label{appendix:self_consistency_checks}

Algorithms \ref{alg:self_consistency_optimal_wealth} and
\ref{alg:self_consistency_optimal_control} below are the fully-detailed
implementations of the self-consistency checks that are sketched in
Algorithm \ref{alg:self_consistency_checks}.  In both cases we

\begin{enumerate}
  \item simulate $M$ independent paths of the dual process
    $\{Y_{t_n}\}_{n=0}^N$ via \eqref{eq: dual process},
  \item stop each path at the first time
    \[
      \tau^* \;=\;\min\{\,t_n:\,\widetilde U_K(Y_{t_n})\ge\widetilde V(t_n,Y_{t_n})\}\wedge T,
    \]
  \item compute the discounted utility
    \[
      P_i \;=\;e^{-\beta\,\tau^{*(i)}}\,U\bigl(X_{\tau^{*(i)}}^{(i)}-K\bigr),
      \qquad
      r_i \;=\;\frac{P_i - V_0}{V_0},
    \]
  \item and finally report the sample mean $\bar r$ and standard deviation
    $s_r$ of $\{r_i\}_{i=1}^M$.
\end{enumerate}

The only difference between the two algorithms is in how the terminal
primal wealth $X_{\tau^*}$ is obtained:

\begin{itemize}
  \item In Algorithm \ref{alg:self_consistency_optimal_wealth}, it is recovered analytically by $X_{\tau^*} =-\,\partial_y\widetilde V\bigl(\tau^*,Y_{\tau^*}\bigr)$
    as in \eqref{eq: optimal wealth}.
  \item In Algorithm \ref{alg:self_consistency_optimal_control}, it is evolved forward under the control $\pi_t^*=\tfrac{\theta}{\sigma}\,Y_t\,\widetilde V_{yy}(t,Y_t)$
    as in \eqref{eq: optimal control}.
\end{itemize}

The reference value
\[
  V_0 \;=\;V(0,x_0)
      \;=\;\min_{y\in(e^{-C_z},\,e^{C_z})}\Bigl\{\widetilde V(0,y)+x_0\,y\Bigr\}
\]
is exactly the dual-primal transform introduced in Algorithm \ref{alg:primal_from_network} and Section \ref{sec:numerical 1}.
From Theory \ref{thm: verification theorem}, the three approaches should give the same value $V_0$ for the same initial wealth $x_0$.

We use 5 seeds for the Monte-Carlo simulation in the self-consistency checks. The seeds are 60, 61, 62, 63, and 64.

\begin{algorithm}[htbp]
\caption{Self-consistency check via optimal wealth from \eqref{eq: optimal wealth} (cf.\ Section \ref{sec:numerical 2})}
\label{alg:self_consistency_optimal_wealth}
\begin{algorithmic}[1]
\Require Trained network $\widetilde{V}(t,y)$; model parameters 
  $x_0,r,\beta,\sigma,\mu,T,\theta,K,\gamma$; utility $U$; 
  time grid $\{t_n = n\,\Delta t\}_{n=0}^N$ with $N=T/\Delta t$; paths $M$.  
\State Precompute $V_0 \!\leftarrow\! V(0,x_0)$.
\For{$i=1,\dots,M$}
  \State {\bfseries (Step 1)} 
    Compute initial dual
    \[y_0^{(i)} \;\leftarrow\;\arg\min_{y>0}\{\widetilde V(0,y)+x_0\,y\}.\]
  \State {\bfseries (Step 2)} 
    \For{$n=1,\dots,N$}
      \State Simulate $Z\sim N(0,1)$ and set
      \[Y_{t_n}^{(i)} \;\leftarrow\; Y_{t_{n-1}}^{(i)}\,
         \exp\bigl((\beta-r-\tfrac12\theta^2)\Delta t
         \;-\;\theta\sqrt{\Delta t}\,Z\bigr).\]
      \If{$\widetilde V(t_n,Y_{t_n}^{(i)}) \le \widetilde U_K(Y_{t_n}^{(i)})
             \;\lor\; n=N$}
        \State $n^*\leftarrow n$; \textbf{break}
      \EndIf
    \EndFor
    \State $\tau^{*(i)}\leftarrow t_{n^*}$.
  \State {\bfseries (Step 3)} 
    Compute stopped wealth
    \[X^{(i)}_{\tau^{*(i)}} \;\leftarrow\;
      -\,\partial_y\widetilde V\bigl(\tau^{*(i)},Y_{\tau^{*(i)}}^{(i)}\bigr).\]
  \State {\bfseries (Step 4)} 
    \[P_i \;\leftarrow\; e^{-\beta\,\tau^{*(i)}}\,
       U\!\bigl(X^{(i)}_{\tau^{*(i)}}-K\bigr),\qquad
       r_i \;\leftarrow\;\frac{P_i - V_0}{V_0}.\]
\EndFor
\State {\bfseries (Step 5)} 
  \[\bar r \;\leftarrow\;\frac1M\sum_{i=1}^M r_i,\qquad
    s_r \;\leftarrow\;\sqrt{\frac{1}{M-1}\sum_{i=1}^M(r_i-\bar r)^2}.\]
\State \textbf{Output:} Sample mean $\bar r$ and standard deviation $s_r$ of the relative errors.
\end{algorithmic}
\end{algorithm}

When simulating the wealth process $X_t$ via its dynamics
\[
dX_t = r X_t \, dt + \sigma \pi_t \left( \theta \, dt + dW_t \right),
\]
we will use the following log-transformed wealth process to avoid $X_t \leq K$. Define $\widetilde{X}_t=\ln(X_t-K)$, then
\[
d\widetilde{X}_t=\frac{rX_t\,dt+\sigma\pi_t\left(\theta\,dt+dW_t\right)}{X_t-K}-\frac{\sigma^2\pi_t^2\,dt}{2(X_t-K)^2}.
\]
$X_t$ is then replaced by $K+\exp(\widetilde{X}_t)$ in the algorithm.

\begin{algorithm}[htbp]
\caption{Self-consistency check via optimal control from \eqref{eq: optimal control} (cf.\ Section \ref{sec:numerical 3})}
\label{alg:self_consistency_optimal_control}
\begin{algorithmic}[1]
\Require Trained network $\widetilde V(t,y)$; model parameters 
  $x_0,r,\beta,\sigma,\mu,T,\theta,K,\gamma$; utility $U$; 
  time grid $\{t_n = n\,\Delta t\}_{n=0}^N$ with $N=T/\Delta t$; paths $M$.
\State Precompute $V_0 \!\leftarrow\! V(0,x_0)$.
\For{$i=1,\dots,M$}
  \State {\bfseries (Step 1)} 
    \[y_0^{(i)} \leftarrow \arg\min_{y>0}\{\widetilde V(0,y)+x_0\,y\},\quad
      X_0^{(i)} \leftarrow x_0.\]
  \State {\bfseries (Step 2)} 
    \For{$n=1,\dots,N$}
      \State Simulate $Z_n\sim N(0,1)$ and set
      \[Y_{t_n}^{(i)} = Y_{t_{n-1}}^{(i)}
         \exp\bigl((\beta-r-\tfrac12\theta^2)\Delta t
         - \theta\sqrt{\Delta t}\,Z_n\bigr).\]
      \If{$\widetilde V(t_n,Y_{t_n}^{(i)}) \le \widetilde U_K(Y_{t_n}^{(i)})
            \;\lor\; n=N$}
        \State $n^*\leftarrow n$; \textbf{break}
      \Else
        \State Compute control
        \[\pi_n^{*(i)} \;\leftarrow\; \frac{\theta}{\sigma}\,Y_{t_n}^{(i)}\,
           \widetilde V_{yy}\bigl(t_n,Y_{t_n}^{(i)}\bigr).\]
        \State Update log-transformed wealth
        \[
          \widetilde X_{n}^{(i)}
          = \ln\bigl(X_{n-1}^{(i)}-K\bigr)
            +\frac{rX_{n-1}^{(i)}\Delta t
              +\sigma\,\pi_n^{*(i)}(\theta\,\Delta t+Z_n\sqrt{\Delta t})}
                 {X_{n-1}^{(i)}-K}
            -\frac{\sigma^2(\pi_n^{*(i)})^2\Delta t}
                   {2\,(X_{n-1}^{(i)}-K)^2}.
        \]
        \State Update wealth
        \[
          X_{n}^{(i)} \leftarrow K + \exp\bigl(\widetilde X_{n}^{(i)}\bigr).
        \]
      \EndIf
    \EndFor
    \State $\tau^{*(i)}\leftarrow t_{n^*}$, \quad $X^{(i)}_{\tau^{*(i)}}\leftarrow X_{n^*}^{(i)}$.
  \State {\bfseries (Step 3)} 
    \[P_i \leftarrow e^{-\beta\,\tau^{*(i)}}\,U\bigl(X^{(i)}_{\tau^{*(i)}}-K\bigr),
      \quad r_i \leftarrow \frac{P_i - V_0}{V_0}.\]
\EndFor
\State {\bfseries (Step 4)} 
  \[\bar r \leftarrow \frac1M\sum_{i=1}^M r_i,\qquad
    s_r \leftarrow \sqrt{\frac{1}{M-1}\sum_{i=1}^M(r_i-\bar r)^2}.\]
\State \textbf{Output:} Sample mean $\bar r$ and standard deviation $s_r$ of the relative errors.
\end{algorithmic}
\end{algorithm}

\end{document}